\documentstyle[11pt]{article}
\font\mtbb=msbm10 scaled \magstep1
\def\R{\hbox{\mtbb R}}
\def\C{\hbox{\mtbb C}}
\def\Z{\hbox{\mtbb Z}}
\def\CA{{\cal A}}

\def\bt{\bullet}
\def\pt{\partial}

\catcode`@=11
\def\ra#1{\;\; \vbox{\m@th\ialign{##\crcr
     \hfil ${}^{\scriptstyle #1}$\hfil\crcr\noalign{\kern-\p@\nointerlineskip}
      $\hfil\rightarrow\hfil$\crcr}} \;\;}

\def\ri#1{\;\; \vbox{\m@th\ialign{##\crcr
     \hfil ${}^{\scriptstyle #1}$\hfil\crcr\noalign{\kern-\p@\nointerlineskip}
      $\hfil\hookrightarrow\hfil$\crcr}} \;\;}

\catcode`@=12 

\def\<{\langle}
\def\>{\rangle}
\def\nn{\nonumber}

\def\id{\hbox{id}}

\def\beq{\begin{equation}}
\def\eeq{\end{equation}}
\def\bea{\begin{eqnarray}}
\def\eea{\end{eqnarray}}

\def\l#1{\label{#1}}

\def\={\; = \;}

\def\id{\hbox{id}}
\def\kr{\hbox{Ker\ }}

\def\pt{\partial}
\def\bt{\wedge}
\def\={\; = \;}

\def\bo{\begin{observation}}
\def\eo{\end{observation}}
\def\bxx{\begin{example}}
\def\exx{\end{example}}

\begin{document}
\thispagestyle{empty}
\null
\vspace{5cm}

\def\thefootnote{\fnsymbol{footnote}}
\begin{center}
{\LARGE Field Theory on $q=-1$ Quantum Plane}

\ \\
\ \\

Andrzej Sitarz \footnote{Partially supported
by KBN grant 2P 302 103 06 } \\
\ \\
{\em Department of Theoretical Physics \\
Institute of Physics, Jagiellonian University \\
Reymonta 4, 30-059 Krak\'ow, Poland \\
e-mail: sitarz@if.uj.edu.pl}
 \\

\end{center}

\vfill
\begin{abstract}

We build the $q=-1$ defomation of plane on a product
of two copies of algebras of functions on the plane.
This algebra constains a subalgebra of functions on
the plane. We present general scheme (which could be
used as well to construct quaternion from pairs of complex
numbers) and we use it to derive differential structures,
metric and discuss sample field theoretical models.
\end{abstract}

\vfill
\noindent{\sc TPJU 4/95}\\
{\sc February 1995 }

\vfill
\newpage

\parindent 0pt
\def\htt#1{\widehat{#1}}

\section{Introduction}

The algebraic concept of the two-dimensional plane,
conceived as the $C^*$ algebra of continuous functions
it, can be easily generalised to the noncommutative case.
The $q$-deformation of such algebra, called Manin plane,
is the simplest example (see \cite{MAN} for basic introduction,
also \cite{CON} for general concepts of noncommutative geometry).
Let us remind here its definition, we construct this algebra from the
self-adjoint generators $X$ and $Y$ and the unit $1$,
satisfying the following commutation relations:
\begin{equation}
XY \; = \; q YX,
\label{qq1} \end{equation}
where $q$ is a non-zero complex number. Since we can choose the
generators to be selfadjoint we have $q q^*=1$.

In this paper we should restrict ourselves to one
particular case when $q=-1$, and we shall call the
 algebra in this situation as corresponding to
'noncommutative plane' \cite{ZAK}.

The quantum planes in general have been the topic of several
papers \cite{WES,ZUM}, focusing especially on the differential
structures. What we shall discuss here is a special formulation of
$q=-1$ quantum plane algebra, which is built on the product
of algebras of functions on the standard plane. In this respect our
approach differa from \cite{ZAK}, moreover, we do not discuss the
corresponding $q=-1$ quantum linear group. We shall
describe the general construction, demonstrating that it
could be applied also to some other objects. We shall give
simple examples of this scheme (obtaining in this way quaternions
{}from complex numbers, for instance).

The main body of this paper is devoted, however, to the application
of this method to $q=-1$ deformation of the plane and simple models
of field theory in this setup. We discuss the differential calculus,
metric and construct sample actions for scalar field theory and
gauge theory.

\section{General concepts}

In this section we shall introduce a general structure of
an algebra built on the product of two copies of an
arbitrary algebra. We shall see that this structure could
be used to constrict the algebra of quaternions as well as
the $q=-1$ deformation of the algebra of functions on
the plane.

Let $\CA$ be an algebra and $\htt{\ } : \CA \to \CA$ be
and algebra automorphism such that $\htt{\ } \circ \htt{\ } = \id$.

Then the following is true:
\bo{1}
If $\xi$ is an element of the centre of $\CA$, which is $\htt{\ }$ invariant,
$\htt{\xi} = \xi$, then $\CA \times \CA$ with pointwise addition and
the following multiplication rule:
\begin{equation}
(a,b) \bt (A,B) \= (a A + \xi b \htt{B} , b \htt{A} + a B)
\end{equation}
is an algebra.
\eo
As linearity is quite obvious, it remains to check the associativity,
the proof of which is technical. We shall call this algebra $htt{\CA}$.

Additionally, if $\CA$ has a $\star$-structure and $\xi^\star=\xi$
we can choose between the following star structures on $\htt{\CA}$:
\begin{equation}
(a,b)^\star \= (a^\star, \htt{b}^\star), \l{fss}
\end{equation}
or:
\begin{equation}
(a,b)^\star \= (a^\star, -\htt{b}^\star). \l{sss}
\end{equation}

\begin{example}
Let $\CA$ be just the field of complex numbers. Then if \ $\htt{}$\
is just the complex conjugation (as $\C$ is commutative then there is
no problem with it being an algebra automorphism) and $q=-1$, then
$\htt{\CA}$ with the second $\star$-structure (\ref{sss})
is the algebra of quaternions.
\end{example}

\bo{2} $\CA$ is included in $\htt{\CA}$:
\beq \CA \ni a \; \mapsto \; (a,0) \in \htt{\CA}, \eeq
\eo

To end this section let us present two examples of construction
based on the algebras of functions on discrete spaces. In the
first case we take as $\CA$ the algebra of complex
valued functions on $\Z_2$, and $\Z_2 \times \Z_2$ in the second one.

The algebra of functions on $\Z_2$ has one generator
$a$ such that $a^2=1$, $a^\star=a$. Let us take $\htt{a} = -a$
(this operation corresponds to underlying symmetry of changing
the points), then we must choose $\xi$ to be real (without loss
of generality it could be set to $1$) and we get the algebra
$\htt{\CA}$ generated by $a=(a,0)$ and $b=(0,1)$ with the relations:
\beq
a^2 = 1, \;\;\;\; b^2=1 \;\;\;\; ab + ba = 0. \l{al-1}
\eeq

Now, let us turn to $\Z_2 \times \Z_2$, functions on this space
are generated by commuting $a$ and $b$, $a^2=b^2=1$. Let us
take $\htt{a}=a$, $\htt{b}=-b$ and $\xi=a$. Then the algebra
$\htt{\CA}$ is generated by $a=(a,0)$ and $B = (0,1)$ with the
following relations:
\beq
b^2=1, \;\;\;\; A^2 = a, \;\;\;\; Ab + bA = 0. \l{al-2}
\eeq
Next we shall apply this procedure to obtain $q=-1$
Manin plane.

\section{Anticommutative plane}

Let us take $\CA$ to be an algebra of polynomials on $\R^2$ and
let us denote its generators by $x$ and $y$. We choose \ $\htt{}$\ \
to be the automorphism of mirror symmetry: $\htt{x}=x$, $\htt{y}=-y$
and, finally, $\xi=x$. Then the algebra $\htt{\CA}$ has two generators:
$Y = (y,0)$ and $X=(0,1)$ with the following relations:
\beq XY + YX \= 0. \l{qplane} \eeq
Of course, we have $X^2=x$. The relation (\ref{qplane}) is just the deformation
relation of the Manin plane for $q=-1$. What have we achieved here is that now
we
can work with the presentation of such algebra as an algebra built
upon the cartesian product of algebras of functions on the regular
plane.

The choice of the $\star$-structure corresponds now to the choice
of $X$ to be a hermitian or antihermitian operator (of, course in
either case $X^2=x$ would be hermitian.

\subsection{Differential calculus}

In this section we shall construct a differential calculus
on our anticommutative plane. From many possibilities we
shall choose the one, which has the following properties.
First, the bimodule of one-forms must be a free module
generated by $dX$ and $dY$, moreover, we require that
the calculus, when restricted to $\CA \subset \htt{\CA}$,
should be just a standard differential calculus on the plane.

This gives us immediately:
\bea
 Y\; dY = dY\; Y & \;\;\;\;\;\;\; & x\; dY = dY\; x \l{c1} \\
 x\; dx = dx\; x & \;\;\;\;\;\;\; & Y\; dx = dx\; Y \l{c2}
\eea
where $dx = d(X^2) = X\;dX+dX\;X$.

If we look for calculi, which are defined through relations
of the type $x^k dx^l = C_{ij}^{kl} dx^i x^j$ and assume
(\ref{c1}-\ref{c2}) we obtain that there are only three
possibilities of differential structures:
\labelsep 0.5cm
\begin{center}
\begin{itemize}
\item[(A)]
$X\; dX = dX\; X,$\\
$Y\; dX = - dX\; Y,$\\
$X\; dY = - dY\; X$
\item[(B)]
$X\; dX = - dX\; X + w dY\; Y, $ \\
$Y\; dX = - dX\; Y,$\\
$X\; dY = - dY\; X,$
\item[(C)]
$X\; dX =  dX\; X + w dY\; Y,$\\
$X\; dY = dY\; X,$\\
$Y\; dX = - dX\; Y - 2dY\; X.$
\end{itemize}
\end{center}
The first one (A) is the standard one, as
it is built from two commutative calculi (one with
$X\, dX=dX\, X$, the other one with $Y\,dY=dY\,Y$).
This is also the only calculus out these three that
has $\kr d = \C$. Indeed, let us observe that for
the other calculi from the relations (B-C)
we immediately get:
\beq d ( X^2 - \frac{w}{2} Y^2 ) \= 0. \eeq

In either of the presented versions of differential calculi
we could have the higher order calculus:
\beq
dX \bt dX =  0, \;\;\;\;\;\;\;\;\;
dY \bt dY  =  0,  \;\;\;\;\;\;\;\;\;
dX \bt dY  =   dY \bt dX, \l{calch}
\eeq
however, only in the first one (A) the condition
$dX \wedge dX=0$ is necessary.

In our further considerations we shall use only the first
calculus, (however, we shall briefly discuss the existence
of the metric structures on all of them), mainly due to the
fact that it is the only one with $\kr d = \C$, as we would
normally assume for physical theories.

\subsection{Metric}

We use here the definition of metric as proposed and
discussed in the general case of quantum plane \cite{JA}.
Let us denote $g_{ab}$ the value of metric $g(da,db)$,
$a,b$ being $X$ and $Y$.

It appears that only for the first two calculi (A,B) we can
have a nonvanishing metric, for the third
one, the bimodule properties of the metric enforce that it
should vanish.

We find that for the calculi (A) and (B)
$g_{XX}$ and $g_{YY}$ must be in the centre of $\htt{\CA}$ whereas
$g_{xy}$ and $g_{yx}$ must be of the form $xy f$, where
$f$ is in the centre of the algebra. Additionally, for (B)
we must require that $g_{XY} = - g_{YX}$  and $w g_{yy} = 2 X^2 f$.

\section{Scalar Field Theory}

To construct a simple scalar field theory on $q=-1$ quantum plane,
with the field $\Phi \in \htt{A}$ we need to choose the differential
calculus as well as the $\star$ operation. As we have already dealt with
the first problem, having chosen the calculus (A), we shall show
the result for each $\star$ structure, using the notation: $(a,b)^\star =
(a^\star, \pm \htt{b}^\star)$.

First, for $\Phi$ presented as $(\psi, \phi)$ according to the general
scheme (see Observation 1.) and the calculus (A) we get:
\beq
d (\psi, \phi) \= dX\; (\htt{\phi} + 2 x \pt_x \htt{\phi}, 2 \pt_x \htt{\psi})
+
dY\; (\pt_y \psi, \pt_y \phi), \l{dphi1}
\eeq

Let us mention here that, formally, in the above equations,
$d (X \frac{1}{\sqrt{x}}) = 0$, however, this does not pose
a problem as $\frac{1}{\sqrt{x}}$ does not belong to the algebra
of functions, which we are considering.

Having calculated $d \Phi$ we can for any metric $g$ calculate
$g(d \Phi^\star, \Phi)$, which is the kinetic term of the
Lagrange function of any field theory. We shall do it for the
simplest admissible metric $g_{XX}=g_{YY}=1$ and $g_{XY}=g_{YX}=0$,
we use also notation $2 \pt_X f = f + 2 \pt_x f$. Then for (\ref{dphi1}) one
obtains:
\bea
g(d \Phi^\star,d\Phi) & = & 4 (\pt_X \htt{\phi}^\star, \pm \pt_x \psi^\star)
(\pt_X \htt{\phi},  \pt_x \htt{\psi}) \nn \\
 & + & (\pt_y \psi^\star, \pm \pt_y \htt{\phi}^\star) (\pt_y \psi, \pt_y \phi)
\eea

which leads to the following final expression:

\bea
&& g(d \Phi^\star, d\Phi)  =  4|\pt_X \htt{\phi}|^2 \pm 4 x |\pt_x \psi|^2
+ | \pt_y \psi|^2 \pm x |\pt \htt{\phi} |^2 \nn \\
&& \;\;\; + \left( 4 ( \pt_X \htt{\phi}^\star \pt_x \htt{\psi} \pm
\pt_x {\psi}^\star \pt_X \phi ) + ( \pt_y \psi^\star \pt_y \phi
\pm \pt_y \htt{\phi}^\star \pt_y \htt{\psi}) \right) X. \l{meme}
\eea

For fields $\Phi$, which belong to the subalgebra of functions on
the regular plane $\Phi_r = (\psi,0)$ the above expression reduces to:

\beq
g(d \Phi_r^\star, d\Phi_r) = \pm 4 x |\pt_x \psi|^2 + |\pt_y \psi|^2,
\eeq

so it is the standard kinetic term on a plane with a nonconstant metric.

For fields, which are of the form $\Phi_n = (0,\phi)$, the kinetic
term (\ref{meme}) becomes:
\beq
g(d \Phi_n^\star, d\Phi_n) =  4 x^2 |\pt_x \htt{\phi}|^2 +
x |\pt_y \htt{\phi}|^2 + |\htt{\phi}|^2 +
2x (\htt{\phi}^\star \pt_x \htt{\phi} + \htt{\phi} \pt_x \htt{\phi}^\star).
\eeq

\subsection{Integration and Action}

The concept of the integration (or trace) on the algebra is one of crucial
elements, which are required to build a satisfactory model. Though it is rather
unclear in general whether such operation could be introduced for
quantum spaces, even if it could be a number-valued operation, in this
particular case of $q=-1$ deformation one formulate the problem without
difficulties. We propose the trace to be a $\C$-valued, symmetric operator
on our algebra $\htt{\CA}$, which is $\star$-covariant, i.e.
$\int (ab) = \int(ba)$ and $\int a^\star = (\int a)^\star$. Suppose we have
a trace $\int_0$ on the algebra $\CA$ (remember that $\htt{\CA}$ is built on
the product of two copies of $\CA$), which is $\ \htt{} \ $ invariant. Then
the following is true:

\bo{3} $ \htt{\CA} \ni (\psi,\phi ) \to \int_0 \psi $ is a trace on the
algebra $\htt{\CA}$.
\eo

For the algebra of quaternions, as discussed at the beginning of our paper,
this gives the natural answer: the $\star$ invariant trace on the algebra
of complex numbers is the real part: $ \int_0 z  = \hbox{Re\ } z $,
consequently,
a trace of the quaternion $(a,b)$ is just the real part of $a$.

\medskip

Coming back to the algebra of $q=-1$ quantum plane we can now write the
proposition for the simplest action of scalar field theory as the trace
of (\ref{meme}):
\bea
{\cal S} & = & \int_0 \left(
4 x^2 |\pt_x \htt{\phi}|^2
+ x (\htt{\phi}^\star \pt_x \htt{\phi} + \htt{\phi} \pt_x
\htt{\phi}^\star) \right). \nn \\
& \pm & \left(. 4 x |\pt_x \psi|^2 + | \pt_y \psi|^2
\pm x |\pt_y \htt{\phi} |^2 + |\htt{\phi}|^2 \right)
\eea
where $\int_0$ is any integration on the plane, which is $\ \htt{} \ $
invariant. As we can see, we have this action could be interpreted as
a standard action (though in a non-trivial metric) for a doublet of scalar
fields on the plane, with the usual kinetic terms as well as a mass term
for one field.

\section{Gauge Theory}

In this section we shall briefly discuss the construction of a simple
gauge theory on $q=-1$ noncommutative plane. We shall use the concepts
of Connes' approach \cite{CON} to gauge theories in noncommutative geometry,
choosing the gauge group as the unitary group of the algebra, being aware
that this scheme may not be appropriate for arbitrary quantum deformation.

The unitary group ${\cal U}(\htt{\CA})$ consist of all elements $(u,v)$ such
that:
\beq
u u^\star \pm x v v^\star = 1, \;\;\;\;\;\;\;\;
v \htt{u}^\star \pm u \htt{v}^\star = 0.
\eeq

Of course, the group all unitary functions on the plane is a subgroup
of ${\cal U}(\htt{\CA})$.

The connection one-form $A$ must be antihermitian $A = -A^\star$, so
that the curvature $F = dA + A \wedge A$ is hermitian. The two-form
$F$ has, due to the rules of calculus (\ref{calch}), has only one
component $dX \wedge dy F_{xy}$.

If $A = dX ( \psi_X, \phi_X) + dy (\psi_y, \phi_y)$ the relation
$A^\star = -A$ gives us:
\bea
\psi_X = - \htt{\psi}_X^\star & \;\;\;\; &
\phi_X = \mp \phi_X^\star \l{conjcond1} \\
\psi_y = - \psi_y^\star & \;\;\;\; &
\phi_y = \pm \htt{\phi}_y^\star \l{conjcond2}
\eea
and $F$ becomes:
\bea
F & =  &  dX \wedge dY \left( - \pt_y \psi_X
 \pm \phi_y^\star \pm 2x \pt_x \phi_y^\star
\psi_X ( \psi_y + \htt{\psi}_y ) + x \phi_y^\star
(\htt{\phi_X}  - \phi_X) \right) \nn \\
& & + dX \wedge dY \left( - \pt_y \phi_X + 2 \pt_x \psi_y^\star
+ \psi_X \phi_y \mp \psi_X^\star \phi_y^\star \right) X
\eea

The expression for $F_{xy}^\star F_{xy}$ and its integral (which
would be the Yang-Mills action) is so complicated that it is rather
unreadable and we shall not discuss it in detail. Let us observe, however,
that it would certainly contain a potential term for fields $\phi$ and $\psi$,
which is a feature that distinguishes such theory from its analogue on the
commutative plane. The potential term $V$ would be:
\beq V \; \sim \;
x \left( \psi_X \phi_y \mp \psi_X^\star \phi_y^\star \right)^2.
\eeq

To end this section we shall briefly present the calculation of
$F$ in situation, where we restrict ourselves to such connections
that arise from the usual $U(1)$ gauge group, which have the
following form:
\beq A_r = dX ( 0, \phi_X) + dy (\psi_y, 0) \eeq
where (\ref{conjcond1},\ref{conjcond2}) for $\phi_X$ and $\psi_y$
still hold. The curvature now reads:
\beq
F  \=   dX \wedge dy \left(
- \pt_y \phi_X + 2 \pt_x \psi_y^\star \right)
\eeq
and its square gives (up to the factor 2) the usual action for the
two-dimension abelian Yang-Mills theory.

\section{Conclusions}

As we have demonstrated one could construct the $q=-1$ deformation
in such a way that it is just an extension of the algebra of functions
on the standard plane in the same way as quaternions are extension of
complex numbers. It is unclear, whether such property holds only for
$q=-1$ or could be generalised for other $q$ (roots of unity in
particular).

Another interesting problem is the application the scheme to some other
algebras, which could lead to the generalisation of the concept of $q=-1$
deformations for other objects (including discrete spaces).

Discussing the examples of field theory models on the anticomutative plane,
we have found that they give some effective actions for usual (commutative)
fields on the plane in a nontrivial metric background. The interesting
features
are the mass terms appearing automatically for scalar field action and
the potential
term in the gauge theory. Though the physical significance of such
models might not be great, they remain interesting as a testing ground
for implications and role of $q$-deformations in theoretical physics.

\def\v#1{{ \bf  #1} }

\end{document}